\documentclass[a4paper,12pt]{article}
\usepackage[top=2.5cm,bottom=2.5cm,left=2.5cm,right=2.5cm]{geometry}
\usepackage{amsfonts}
\usepackage{latexsym,amsmath,amssymb,amsfonts,epsfig,psfrag,url,graphics,ifpdf,multicol}
\usepackage{color,colordvi,amscd,amsthm} %eepic
\usepackage{tikz}
\usepackage[numbers,sort&compress]{natbib}
\usepackage{indentfirst,graphics,epsfig}
\usepackage{graphicx}
\usepackage{graphics}
\usepackage{graphicx,psfrag}
\usepackage{graphics}
\usepackage{mathrsfs}
\usepackage{tabularx}
\usepackage{booktabs}
\usepackage{floatrow}
\usepackage{multirow}
\usepackage{graphicx}
\usepackage{caption}
\usepackage[ruled]{algorithm2e}
\newfloatcommand{capbtabbox}{table}[][\FBwidth]

\newtheorem{theo}{Theorem}[section]
\newtheorem{lem}[theo]{Lemma}
\newtheorem{cl}{Claim}
\newtheorem{coro}[theo]{Corollary}

 \theoremstyle{definition}

\theoremstyle{remark}

\def\pf{\noindent {\bf Proof.} }

\allowdisplaybreaks
\usepackage{indentfirst,subfig}
\begin{document}
	%\thispagestyle{empty}
	%\rule{0cm}{0.5mm}
	\captionsetup[figure]{labelfont={bf},name={Fig.},labelsep=period}

\begin{center} {\large
Maximum dissociation sets in subcubic trees}
\end{center}

\begin{center}
{
  {\small Lei Zhang$^a$, Jianhua Tu$^{a,}$\footnote{Corresponding author.\\\indent \  {\it E-mail addresses:} 2018200896@mail.buct.edu.cn (L. Zhang), tujh81@163.com (J. Tu), xinchl@mail.buct.edu.cn (C. Xin).}, Chunlin Xin$^b$}\\[2mm]

{\small $^a$ Department of mathematics, Beijing University of Chemical Technology, \\
\hspace*{1pt} Beijing, P.R. China 100029} \\[2mm]
{\small $^b$ School of Economics and Management, Beijing University of Chemical Technology, \\
\hspace*{1pt} Beijing, P.R. China 100029} \\[2mm]
}
\end{center}

\begin{center}

\begin{abstract}
A subset of vertices in a graph $G$ is called a maximum dissociation set if it induces a subgraph with vertex degree at most 1 and the subset has maximum cardinality. The dissociation number of $G$, denoted by $\psi(G)$, is the cardinality of a maximum dissociation set. A subcubic tree is a tree of maximum degree at most 3. In this paper, we give the lower and upper bounds on the dissociation number in a subcubic tree of order $n$ and show that the number of maximum dissociation sets of a subcubic tree of order $n$ and dissociation number $\psi$ is at most $1.466^{4n-5\psi+2}$.

\vspace{5mm}

\end{abstract}
%\end{minipage}
\end{center}

\baselineskip=0.24in
\section{Introduction}
We consider only finite, simple, and undirected labeled graphs, and use
Bondy and Murty \cite{Bondy2008} for terminology and notations not
defined here.

In a graph $G$, an \emph{independent set} is a set of pairwise non-adjacent vertices of $G$. An independent set is \emph{maximal} if it is not a proper subset of any other independent set, and \emph{maximum} if it has maximum
cardinality. The \emph{independence number} of a graph $G$, denoted
by $\alpha(G)$, is the cardinality of a maximum independent set of
$G$.

A dissociation set in a graph $G$ is a vertex subset $F$ such that the subgraph $G[F]$ induced by $F$ has vertex degree at most 1.
A \emph{maximum dissociation set} of $G$ is a dissociation set of maximum
cardinality. The \emph{dissociation number} of a graph $G$, denoted
by $\psi(G)$, is the cardinality of a maximum dissociation set
of $G$. The problem of finding a maximum dissociation set in a given graph has been introduced by Yannakakis \cite{Yannakakis1981} and is known to be NP-hard for bipartite graphs. The complexity of the problem for some classes of graphs has been studied \cite{Alekseev2007,Cameron2006,Orlovich2011,Xiao2017, Yannakakis1981}. Note that a set $F$ of vertices of a graph $G$ is a dissociation set if and only if its complement $V(G)\setminus F$ is a so-called 3-path vertex cover, that is, a set of vertices of $G$ intersecting every path of order 3 in $G$. The 3-path vertex cover problem is to find a minimum 3-path vertex cover in a graph $G$ and has received considerable attention in the literature \cite{Bresar2011,Kardos2011,Katrenic2016,Tu2011,Xiao2017}.

In 1960s, Erd\H{o}s and Moser raised the problem of determining the maximum number of maximal independent sets for a general graph $G$ of order $n$. This problem was solved by Erd\H{o}s, and later Moon and Moser \cite{Moon1965}. Since then, the problem was extensively studied for various classes of graphs, including trees \cite{Wilf1986,Sagan1988}, connected graphs \cite{Furedi1987,Griggs1988}, bipartite graphs \cite{Liu1993}, unicyclic connected graphs \cite{Koh2008}, connected graphs with at most $r$ cycles \cite{Sagan2006}. A number of authors have also studied the problem of determining the maximum number of maximum independent sets for various classes of graphs \cite{Zito1991,Jou2000,Sagan2006,Mohr2018,Mohr2018-2,Derikvand2014}. On the other hand, the problem of determining the maximum number of minimal (or minimum) dominating sets was also studied in the literature \cite{Alvarado2019, Brod2006, Connolly2016, Fomin2005}.

%the problem of determining the maximum number of maximum independent sets or maximal independent sets in graphs has been extensively studied. In 1986, Wilf \cite{Wilf1986} determined the maximum number of maximal independent sets in a tree and Sagan \cite{Sagan1988} gave a short proof and determined the extremal trees. In 1991, Zito \cite{Zito1991} determined the maximum number of maximum independent sets in a tree of order $n$. The problem has been studied for various classes of graphs, including connected graphs \cite{Jou2000}, bipartite graph \cite{Liu1993}, connected graphs with at most $r$ cycles \cite{Sagan2006}, etc. A subcubic tree is a tree such that every vertex is incident with at most three edges.

A subcubic tree is a tree of maximum degree at most 3. Recently, Mohr et al. \cite{Mohr2018} considered subcubic trees and proved the following.

\begin{theo}\cite{Mohr2018}
If $T$ is a subcubic tree of order $n$ and independence number $\alpha$, then the number of maximum independent sets in $T$ is at most
$1.618^{2n-3\alpha+1}.$
\end{theo}

Inspired by these aforementioned problems, we consider the analogous problem of determining the maximum number of maximum dissociation sets in a graph $G$.
In \cite{Tu2019}, we determined the maximum number of maximum dissociation sets in a tree of order $n$ and characterized the extremal trees. In the present paper, we consider the dissociation number and the number of the maximum dissociation sets of a subcubic tree. In the next section, we will give the lower and upper bounds on the dissociation number of a subcubic tree of order $n$. In Section 3, we will show that the number of maximum dissociation sets of a subcubic tree of order $n$ and dissociation number $\psi$ is at most $1.466^{4n-5\psi+2}$.

\section{The bounds on the dissociation number of a subcubic tree}

Let $G$ be a graph and $v$ be a vertex in $G$. The neighborhood $N_G(v)$ is the set of vertices adjacent to $v$ and the closed neighborhood $N_G[v]$ is $N_G(v)\cup\{v\}$.  For a subset $X$ of vertices, the induced subgraph $G[X]$ is the subgraph of $G$ whose vertex set is $X$ and whose edge set consists of all edges of $G$ which have both ends in $X$. If $U$ is the set of vertices deleted, the resulting subgraph is denoted by $G-U$. If $U=\{v\}$, we write $G-v$ for $G-\{v\}$. A vertex in a tree $T$ is called a leaf if it has degree exactly one. The neighbor of a leaf is called a {\it support vertex} of $T$. If a support vertex of $T$ is adjacent to at least two leaves, then it is a {\it strong support vertex} of $T$.

%An {\it endvertex} of $G$ is a vertex of degree at most 1. A neighbor of an endvertex of $G$ is a {\it support vertex} of $G$. If a support vertex of $G$ is adjacent to at least two endvertices, then it is a {\it strong support vertex} of $G$.

Let $k$ be a positive integer. A set $S$ of vertices in $G$ is called a $k$-path vertex cover if every path of order $k$ in $G$ contains at least one vertex from $S$. The $k$-path vertex cover number $\tau_k(G)$ of $G$ is the cardinality of a minimum $k$-path vertex cover in $G$.

\begin{theo}{\label{th2.2}}\cite{Bresar2011}
Let $T$ be a tree of order $n$ and $k$ a positive integer. Then, $\tau_k(G)\leq n/k$.
\end{theo}

A set $S$ of vertices of a graph $G$ is a 3-path vertex cover if and only if its complement $V(G)\setminus S$ is a dissociation set. Therefore, we can obtain a lower bound on the dissociation number of a tree of order $n$.

\begin{theo}{\label{th2.2}}
If $T$ is a tree of order $n$, then $\psi(T)\ge \frac{2n}{3}$.
\end{theo}

For a positive integer $\ell$, let $T_{\ell}$ arise by attaching a pendant edge to every vertex of a path of order $\ell$. Hence, $T_{\ell}$ is a subcubic tree of order $3\ell$ and $\psi(T_{\ell})=2\ell$. It follows that for a subcubic tree the lower bound on the dissociation number given by Theorem \ref{th2.2} is tight.

Next, we give the upper bound on the dissociation number of a subcubic tree of order $n$. If $T$ is a tree,
then $T'$ arises from $T$ by attaching a $P_5$ such that $V(T')$ is the disjoint union of $V(T)$ and $\{x,y,z,j,k\}$, and $E(T')=E(T)\cup \{uz,zy,zj,yx,jk\}$, where $u$ is some vertex of $T$.

\begin{theo}\label{th2.3}
If $T$ is a subcubic tree of order $n$, then,
\begin{equation}
\psi(T)\le \frac{4n+2}{5} .\label{1}
\end{equation}
Furthermore, equality holds in (\ref{1}) if and only if $T$ arises from $K_2$ by iteratively attaching $P_5$s.
\end{theo}

\pf Suppose, for a contradiction, that the theorem is false, and let $n$ be the smallest order for which it fails. Let $T$ be a subcubic tree of order $n$ and dissociation number $\psi$ such that either $\psi> \frac{4n+2}{5}$ or $\psi=\frac{4n+2}{5}$ and $T$ doesn't arise from $K_2$ by iteratively attaching $P_5$s.
It is easy to see that $T$ has diameter at least 3. We root $T$ at an endvertex of a longest path in $T$. Let $u$ be a leaf of maximum depth in $T$ and $uvwx$ be a path of $T$.

%let $v$ be the parent of $u$, let $w$ be the parent of $v$ and let $x$ be the parent of $w$.

\setcounter{cl}{0}
\begin{cl}\label{cl2.1}
$d_T(v)=2$.
\end{cl}
\noindent{\bf Proof of Claim \ref{cl2.1}.}
Suppose, for a contradiction, that $d_T(v)=3$. Let $T'=T-(N_T[v]\setminus\{w\})$, then $T'$ has order $n-3$ and dissociation number $\psi-2$. By the choice of $n$, we have
\begin{align}
\psi=\psi(T')+2\le \tfrac{4|V(T')|+2}{5}+2 =\tfrac{4(n-3)+2}{5}+2< \tfrac{4n+2}{5},\notag
\end{align}
which contradicts the choice of $T$.\qed

%Since $\psi(P_3)=2<\tfrac{4\times3+2}{5}$, we assume that $w$ has a parent $x$.
\begin{cl}\label{cl2.2}
$d_T(w)=3$ and $w$ is not a support vertex.
\end{cl}
\noindent{\bf Proof of Claim \ref{cl2.2}.}
If $d_T(w)=2$, then $T'=T-\{u,v,w\}$ has order $n-3$ and dissociation number $\psi-2$. Similarly, we obtain $\psi<\tfrac{4n+2}{5}$, which contradicts the choice of $T$. Therefore, $d_T(w)=3$.

Suppose, for a contradiction, that $w$ is a support vertex. Let $v'$ be the child of $w$ distinct from $v$ and $v'$ is a leaf. Let $T'=T-\{u,v,w,v'\}$, then $T'$ has order $n-4$ and dissociation number $\psi-3$. By the choice of $n$, we obtain
\begin{align}
\psi=\psi(T')+3  \le \tfrac{4|V(T')|+2}{5}+3=\tfrac{4(n-4)+2}{5}+3< \tfrac{4n+2}{5}\notag,
\end{align}
which contradicts the choice of $T$.\qed

By the above two claims, $w$ has a child $v'$ distinct from $v$ and $v'$ has exactly one child $u'$ in $T$. Let $T'=T-\{u,u',v,v',w\}$, then $T'$ has order $n-5$ and dissociation number $\psi-4$. By the choice of $n$, we obtain
\begin{align}
\psi=\psi(T')+4  \le \tfrac{4|V(T')|+2}{5}+4=\tfrac{4(n-5)+2}{5}+4= \tfrac{4n+2}{5}\notag,
\end{align}
which implies that $\psi=\tfrac{4n+2}{5}$ and the tree $T'$ arises from $K_2$ by iteratively attaching $P_5s$. Since $w$ has degree 3 and each of children of $w$ has exactly a child that is a leaf, $T$ also arises from $K_2$ by iteratively attaching $P_5s$. This contradiction completes the proof.\qed

\begin{theo}\label{th2.4}
If $T$ is a subcubic tree of order $n$, then $\frac{2n}{3}\leq \psi(T)\leq \frac{4n+2}{5}$. Moreover, both bounds are tight.
\end{theo}

\section{The number of maximum dissociation sets in a subcubic tree}

In this section, our main result concerns the largest possible value of the number of maximum dissociation set in a subcubic tree of order $n$ and dissociation number $\psi$.

Let $G$ be a graph and $v$ be a vertex of $G$. The set of all maximum dissociation sets of $G$ is denoted by $MD(G)$ and its cardinality by $\Phi(G)$. Let
\begin{align}
\Phi_v(G)&=|\{F\in MD(G): v\in F\}|, \notag\\
\Phi_{\overline{v}}(G)&=|\{F\in MD(G): v\notin F\}|, \notag\\
\Phi_v^0(G)&=|\{F\in MD(G): v\in F \text{ and } d_{G[F]}(v)=0\}|,\notag\\
\Phi_v^1(G)&=|\{F\in MD(G): v\in F \text{ and } d_{G[F]}(v)=1\}|,\notag
\end{align}
It follows that  $\Phi(G)=\Phi_v(G)+\Phi_{\overline{v}}(G)$ and
$\Phi_v(G)=\Phi_v^1(G)+\Phi_v^0(G)$.\par

%Let $f(n)$ denote a sequence, that is
%\begin{center}
%$f(n)=\left \{
%  \begin{array}{ll}
%   \qquad\qquad \qquad1&,  \hbox{if $n=0$;} \\
%   \qquad\qquad\qquad 1&, \hbox{if $n=1$;} \\
%  \qquad\qquad \qquad 3&,  \hbox{if $n=2$;}\\
%    f(n-1)+2f(n-2)+f(n-3)&,\hbox{if $n\ge3$;}
%  \end{array}
%\right.$
%\end{center}
%
%Let $g(n)$ denote another sequence, that is
%\begin{center}
%$g(n)=\left \{
%  \begin{array}{ll}
%   \qquad\qquad \qquad1&,  \hbox{if $n=0$;} \\
%   \qquad\qquad\qquad 2&, \hbox{if $n=1$;} \\
%  \qquad\qquad \qquad 4 &,  \hbox{if $n=2$;}\\
%    g(n-1)+2g(n-2)+g(n-3)&,\hbox{if $n\ge3$;}
%  \end{array}
%\right.$
%\end{center}
%
%By induction, it is easy to get the following two equations.
%\begin{align}
%&g(n)=f(n)+f(n-1)&, \hbox{if $n\ge1$;} \notag  \\
%&f(n)=g(n-1)+g(n-2)&, \hbox{if $n\ge2$.} \notag
%\end{align}\par
\begin{lem}\label{lem3.1}
 Let $T$ be a tree of order $n$ and dissociation number $\psi$. If $v$ is a leaf of $T$ such that $\psi(T-v)=\psi(T)$, then $\Phi_{\overline{v}}(T)\leq \min\{\Phi_v^0(T),\Phi_v^1(T)\}$ and $\Phi_{\overline{v}}(T)\le\frac{1}{3}\Phi(T)$.
\end{lem}
\pf
Let $u$ be the neighbour of $v$ in $T$ and $T'=T-v$. Since $\psi(T')=\psi(T)$, for every maximum dissociation set $F$ of $T'$ the vertex $u$ is contained by the set $F$ and $d_{T'[F]}(u)=1$.
Let $N_{T'}(u)=\{u_1,\cdots,u_k\}$. Next, we will prove that there exists a vertex in $N_{T'}(u)$ such that it is in all maximum dissociation sets of $T'$.

Without loss of generality, suppose to the contrary that there are two maximum dissociation sets $F_1$ and $F_2$ of $T'$ such that $u_1\in F_1$ and $u_2\in F_2$. Let $D_1,\cdots,D_k$ be the connected components of $T'-u$ such that $u_i\in V(D_i)$. See Figure \ref{1}. It is easy to see that $F_1\cap V(D_1)$ is a maximum dissociation set of $D_1$ and every maximum dissociation set of $D_1$ contains $u_1$. Since $F_2\cap V(D_1)$ is a dissociation set of $D_1$ not containing $u_1$, $|F_2\cap V(D_1)|< |F_1\cap V(D_1)|$. Now, $[F_1\cap V(D_1)]\cup [\cup_{i=2}^k(F_2\cap V(D_i))]$
is a maximum dissociation set of $T'$ that does not contain $u$, this is a contradiction. Thus, we verify that there exists a vertex in $N_{T'}(u)$, say $u_1$, such that it is in all maximum dissociation sets of $T'$.

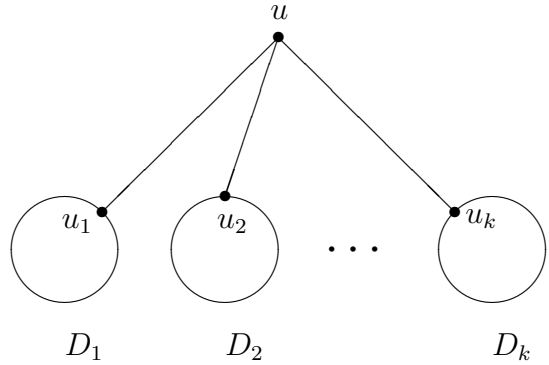
\begin{figure}[h]
	\begin{center}
		\begin{picture}(400,125)
		\multiput(120,35)(60,0){2}{\circle{40}}
        \multiput(280,35)(80,0){1}{\circle{40}}
		\multiput(200,115)(0,0){1}{\circle*{3.5}}
		\multiput(180,55)(0,0){1}{\circle*{3.5}}
		\multiput(266,49)(0,0){1}{\circle*{3.5}}
		\multiput(134,49)(0,0){1}{\circle*{3.5}}
        % \multiput(155,35)(5,0){3}{\circle*{2}}
		\multiput(220,35)(8,0){3}{\circle*{2}}
		\put(200,115){\line(-1,-3){20}}
		\put(200,115){\line(1,-1){66}}
		\put(200,115){\line(-1,-1){66}}
		%%%%%%%%%%%%%%%%%%%%%%%%%%%%%%
		\put(119,43){$u_1$}
		\put(197,122){$u$}
		\put(177,43){$u_2$}
		\put(270,45){$u_k$}
		\put(120,-5){$D_1$}
		\put(180,-5){$D_2$}
		\put(280,-5){$D_k$}
		%%%%%%%%%%%%%%%%%%%%%%%%%%%%
		\caption{\label{2} The tree $T'=T-v$ in the proof of Lemma \ref{lem3.1}}	
		\end{picture}
	\end{center}
\end{figure}

If $F$ is a maximum dissociation set in $T'$, then $F_1=(F\cup\{v\})\setminus\{u\}$ is a maximum dissociation set of $T$ such that $d_{T[F_1]}(v)=0$, and $F_2=(F\cup\{v\})\setminus\{u_1\}$ is a maximum dissociation set in $T$ such that $d_{T[F_2]}(v)=1$. On the other hand, $F$ is also a maximum dissociation set in $T$ that does not contain $v$. Thus, we have
\begin{align}
\Phi(T')=\Phi_{\overline{v}}(T)\leq \min\{\Phi_v^0(T),\Phi_v^1(T)\}\notag.
\end{align}
Since $\Phi(T)=\Phi_{\overline{v}}(T)+\Phi_v^0(T)+\Phi_v^1(T)$, it follows that
$\Phi_{\overline{v}}(T)\le\frac{1}{3}\Phi(T).$ \qed

%Let $n$ and $\psi$ be integers with $\tfrac{2n}{3}\le\psi\le\tfrac{4n+2}{5}$, we denote
%\begin{center}
%$F(n,\psi)=\left \{
%  \begin{array}{ll}
%   1.467^{\frac{4n-5\psi+2}{2}}&,  \hbox{if $\psi$ is even;}\vspace{8pt}\\
%   1.467^{\frac{4n-5\psi+1}{2}}&, \hbox{if $\psi$ is odd.}
%\end{array}
%\right.$
%\end{center}

%Let $T_{(k)}$ be a subcubic tree of order $3k$ defined in the previous section. Since
%\begin{align}
%&\Phi(T_{(1)})=3, \notag\\
%&\Phi(T_{(2)})=6,  \notag\\
%&\Phi(T_{(3)})=13, and  \notag\\
%&\Phi(T_{(k)})=\Phi(T_{(k-1)})+2\Phi(T_{(k-2)})+\Phi(T_{(k-3)})\notag \text{\ for\ every\ } k\ge3,
%\end{align}
%we obtain $\Phi(T_{(k)})=f(k+1)$ for every positive integer $k$.\par

\begin{theo}\label{th3.2}
If $T$ is a subcubic tree of order $n$ and dissociation number $\psi$, then,
\begin{equation}
\Phi(T)\le 1.466^{4n-5\psi+2}.\notag
\end{equation}
\end{theo}

\pf Suppose, for a contradiction, that the theorem is false, and let $n$ be the smallest order for which it fails. Let $T$ be a subcubic tree of order $n$ and dissociation number $\psi$ such that $\Phi(T)>1.466^{4n-5\psi+2}$. It is easy to see that $T$ has diameter at least 3. We can assume that $n\ge4$ and $\psi\ge3$.
%\begin{itemize}
%\item  either $\Phi(T)>F(n,\psi)$
%\item  or $\Phi(T)=F(n,\psi)$  and $\psi=\tfrac{2n}{3}$ but $T \ncong T_{(n-\psi)}$.
%\end{itemize}

%\setcounter{cl}{0}
%\begin{cl}\label{cl3.1}
%The tree $T_0$ contains a path of length at least 3.
%\end{cl}
%\noindent{\bf Proof of Claim \ref{cl3.1}.}
%Suppose, to a contradiction, that $T_0$ is a star $K_{1,n-1}$.\par
%If $n=1$, then $\psi(T_0)=1$ and $\Phi(T_0)=1=g(\frac{4-5+1}{2})$. \par
%If $n=2$, then $\psi(T_0)=2$ and $\Phi(T_0)=1 =f(\frac{8-10+2}{2})$. \par
%If $n=3$, then $\psi(T_0)=2$ and $\Phi(T_0)=3 = f(\frac{12-10+2}{2})$. Furthermore, $T_0=T(1)$.\par
%If $n=4$, then $\psi(T_0)=3$ and $\Phi(T_0)=1< g(\frac{16-15+1}{2})$.\\
%In each case, we obtain a contradiction to the choice of $n$ and $T_0$.\qed

%For simplicity, let $\lambda=1.466$.

Let $\alpha$ be the largest solution of the equation $\alpha^3-\alpha^2-2\alpha-1=0$, that is, $\alpha\approx2.148$. Let $\lambda=\sqrt{\alpha}\approx1.466.$

\setcounter{cl}{0}
%\begin{lem}
% Let $T$ be a subcubic tree of order $n$ and dissociation number $\psi$. If $v$ be a endvertex of $T$ such that $\psi(T-v)=\psi(T)$, then $\Phi(T-v)\le\frac{1}{3}\Phi(T)$.
%\end{lem}
\begin{cl}\label{cl3.1}
Let $T_1$ be a subcubic tree of order $n_1$ and dissociation number $\psi_1$. If $n_1\le n-3$, then for any vertex $v\in V(T_1)$ with $d_{T_1}(v)<3$, we have $\Phi_{\overline{v}}(T_1)\le \lambda^{4n_1-5\psi_1}$.
\end{cl}
\noindent{\bf Proof of Claim \ref{cl3.1}.}
Let $v\in V(T_1)$ and $d_{T_1}(v)<3$. If $v$ is in all maximum dissociation sets of $T_1$, then $\Phi_{\overline{v}}(T_1)=0$. Otherwise, let $T_2$ arise from $T_1$ by attaching a pendant edge $e$ to $v$. Then $T_2$ is a subcubic tree of order $n_1+2$ and dissociation number $\psi_1+2$. If $F$ is a maximum dissociation set of $T_1$ such that $v\notin F$, $F\cup V(e)$ is a maximum dissociation set of $T_2$. By the choice of $n$, this implies
\begin{equation}
\Phi_{\overline{v}}(T_1)\le \Phi(T_2)\le \lambda^{4(n_1+2)-5(\psi_1+2)+2}=\lambda^{4n_1-5\psi_1}.\notag
\end{equation}
We complete the proof of the claim.\qed

We root $T$ at an endvertex of a longest path in $T$. Let $u$ be a leaf of maximum depth in $T$ and $uvwx$  be a path of $T$. For a vertex $z\in V(T)$, let $T_{z}$ be the subtree consisting of $z$ and all of its descendants in $T$.

\begin{cl}\label{cl3.2}
$d_{T}(v)=2.$
\end{cl}
\noindent{\bf Proof of Claim \ref{cl3.2}.}
Suppose, for a contradiction, that $d_{T}(v)=3$. Let $T'=T-V(T_{v})$, then $T'$ has order $n-3$ and dissociation number $\psi-2$. Since $d_{T'}(w)<3$, by Claim \ref{cl3.1}, we have $\Phi_{\overline{w}}(T')\le\lambda^{4(n-3)-5(\psi-2)}=\lambda^{4n-5\psi-2}$. A maximum dissociation set in $T'$ not containing $w$ can be extended in three ways to a maximum dissociation set in $T$, while a maximum dissociation set in $T'$ containing $w$ can only be extended in a unique way to a maximum dissociation set in $T$. Since all maximum dissociation sets in $T$ are of one of these types, we obtain
\begin{align}
\Phi(T)&=3\cdot\Phi_{\overline{w}}(T')+\Phi_w(T')\notag    \\
       &=\Phi(T')+2\cdot\Phi_{\overline{w}}(T')\notag\\
       &\le \lambda^{4(n-3)-5(\psi-2)+2}+2\lambda^{4n-5\psi-2}\notag\\
       &=(\tfrac{1}{\lambda^2}+\tfrac{2}{\lambda^4})\lambda^{4n-5\psi+2}\notag\\
        &<\lambda^{4n-5\psi+2},\notag
\end{align}
where we use $\lambda^2+2<\lambda^4$. This contradiction completes the proof of the claim.\qed

%\begin{cl}\label{cl3.3}
%$w$ is not a support vetrex.
%\end{cl}
%\begin{proof}
%Suppose, for a contradiction, that $w$ is a support vertex. The subcubic tree $T'=T-V(T_v)$ has order $n-2$, while the subcubic tree $T'=T-V(T_w)$ has order $n-4$ and dissociation number $\theta-3$. Now, we consider the following two cases.
%
%\textbf{Case 1.} The subcubic tree $T'$ has dissociation number $\theta-2$.\par
%In this case, every maximum dissociation set in $T$ contains $v$ and $u$. So, we obtain

\begin{cl}\label{cl3.3}
$w$ is not a support vertex.
\end{cl}
\noindent{\bf Proof of Claim \ref{cl3.3}.}
Suppose, for a contradiction, that $w$ is a support vertex. Let $v'$ be a leaf of $T$ that is adjacent to $w$. Let $T^{(1)}=T-V(T_{w})$, then $T^{(1)}$ has order $n-4$ and dissociation number $\psi-3$. Let $T^{(2)}=T-\{u,v\}$, then $T^{(2)}$ has order $n-2$ and dissociation number $\psi-1$ or $\psi-2$. Now, we consider the following two cases.

\textbf{Case 1.} $\psi(T^{(2)})=\psi-1$.\par
In this case, $w$ is in all maximum dissociation sets in $T^{(2)}$. Furthermore, a set $F$ is a maximum dissociation set of $T$ if and only if

\begin{itemize}
\item either $F=F'\cup\{u,v,v'\}$, where $F'$ is a maximum dissociation set of cardinality $\psi-3$ of $T^{(1)}$,
\item or $F=F'\cup\{u\}$, where $F'$ is a maximum dissociation set of cardinality $\psi-1$ of $T^{(2)}$.
\end{itemize}
By the choice of $n$, this implies
\begin{align}
\Phi(T)&= \Phi(T^{(1)})+\Phi(T^{(2)}) \notag \\
         &\le \lambda^{4(n-4)-5(\psi-3)+2}+\lambda^{4(n-2)-5(\psi-1)+2}   \notag\\
         &=(\tfrac{1}{\lambda}+\tfrac{1}{\lambda^3})\lambda^{4n-5\psi+2}      \notag\\
         &=\lambda^{4n-5\psi+2},\notag
\end{align}
where we use $\lambda^2+1=\lambda^3$. This contradiction completes the proof of the case.

\textbf{Case 2.} $\psi(T^{(2)})=\psi-2$.\par
In this case, every maximum dissociation set of $T$ contains both $u$ and $v$. So we have $\Phi(T)=\Phi(T^{(1)})\le\lambda^{4(n-4)-5(\psi-3)+2}<\lambda^{4n-5\psi+2}$. This contradiction completes the proof of the case.

We complete the proof of the claim.\qed

%First we assume $\psi_0$ is even. If $\frac{2n+2}{3}\le\psi_0\le\frac{4n}{5}$, then we have $\Phi(T_0)=\Phi(T')\le g(\frac{4(n-4)-5(\psi_0-3)+1}{2})=g(\frac{4n-5\psi_0}{2})$. By an inductive argument on $n$, it is easy to check $f(n)\ge g(n-1)$ with equality if and only if $n=1$, which contradicts the choices of $T_0$; if $\psi_0=\frac{2n}{3}$, then $\psi_0-3=\frac{2(n-4)-1}{3}$, which contradicts the Theorem \ref{th1}.\par
%
%Next we assume $\psi_0$ is odd, then $\Phi(T_0)=\Phi(T')\le f(\frac{4(n-4)-5(\psi_0-3)+2}{2})=f(\frac{4n-5\psi_0+1}{2})$. By an inductive argument on $n$, it is easy to check $g(n)\ge f(n)$ with equality if and only if $n=0$, which contradicts the choices of $T_0$.\qed

\begin{cl}\label{cl3.4}
$d_{T}(w)=2$
\end{cl}
\noindent{\bf Proof of Claim \ref{cl3.4}.}
Suppose, for a contradiction, that $w$ has a child $v'$ distinct from $v$. By Claim \ref{cl3.2} and \ref{cl3.3}, the vertex $v'$ has exactly one child $u'$ and $u'$ is a leaf. Let $T'=T-V(T_{w})$, then $T'$ has order $n-5$ and dissociation number $\psi-4$. Since every maximum dissociation set of $T$ contains $u$, $u'$, $v$ and $v'$, we have $\Phi(T)=\Phi(T')\le\lambda^{4(n-5)-5(\psi-4)+2}=\lambda^{4n-5\psi+2}$, which contradicts the choice of $T$. This contradiction completes the proof of the claim. \qed

%First we assume $\psi_0$ is even. If $\frac{2n+2}{3}\le\psi_0\le\frac{4n}{5}$, then we have $\Phi(T_0)=\Phi(T')\le f(\frac{4(n-5)-5(\psi_0-4)+2}{2})=f(\frac{4n-5\psi_0+2}{2})$, which contradicts the choice of $T_0$; if $\psi_0=\frac{2n}{3}$, then $\psi_0-4=\frac{2(n-5)-2}{3}$, which contradicts the Theorem \ref{th1}.\par
%
%Next we assume $\psi_0$ is odd. If $\frac{2n+3}{3}\le\psi_0\le\frac{4n+1}{5}$, then we have $\Phi(T_0)=\Phi(T')\le g(\frac{4(n-5)-5(\psi_0-4)+1}{2})=g(\tfrac{4n-5\psi_0+1}{2})$, which contradicts the choices of $T_0$; if $\psi_0=\frac{2n+1}{3}$, then $\psi_0-4=\frac{2(n-5)-1}{3}$, which contradicts the Theorem \ref{th1}.\qed

Since we root $T$ at an endvertex of a longest path in $T$, if the vertex $x$ is the root of $T$, then $T\cong P_4$ and $\Phi(P_4)=2<\lambda^{4\times4-5\times3+2}$. Thus, we can assume that $x$ is a child of a vertex $y$ in $T$.

%Since $\Phi(P_4)=2<\lambda^{4\times4-5\times3+2}$, we may assume that $x$ has parent $y$.

\begin{cl}\label{cl3.5}
$x$ is not a support vertex.
\end{cl}
\noindent{\bf Proof of Claim \ref{cl3.5}.}
Suppose, for a contradiction, that $x$ has a child $w'$ that is a leaf. Let $T^{(1)}=T-V(T_{x})$, then $T^{(1)}$ has order $n-5$ and dissociation number $\psi-3$ or $\psi-4$. Now, we consider the following two cases.

\textbf{Case 1.} $\psi(T^{(1)})=\psi-3$.\par
In this case, every maximum dissociation set of $T^{(1)}$ contains $y$. Let $T^{(2)}=T-\{u,v,w\}$, then $T^{(2)}$ has order $n-3$ and dissociation number $\psi-2$. A set $F$ is a maximum dissociation set of $T$ if and only if
\begin{itemize}
\item either $F\in\{F'\cup\{u,w,w'\},F'\cup\{v,w,w'\}\}$, where $F'$ is a maximum dissociation set of cardinality $\psi-3$ of $T^{(1)}$,
\item or $F=F'\cup\{u,v\}$, where $F'$ is a maximum dissociation set of cardinality $\psi-2$ of $T^{(2)}$.
\end{itemize}
%
%If a maximum dissociation set $F$ of $T$ does not contains $v$, then $F$ must contain both $w$ and $w'$. So there are $\Phi(T')$ maximum dissociation sets in $T$ that contain $u$ but not contain $v$. Similarly, there are $\Phi(T')$ maximum dissociation sets in $T$ that contain $v$ but not contain $u$ and there are $\Phi(T'')$ maximum dissociation sets in $T$ that contain $v$ and $u$. Since all maximum dissociation set in $T_0$ are one of these types,
By the choice of $n$, this implies
%First we assume $\psi_0$ is even. If $\psi_0=\frac{4n-2}{5}$ or $\frac{4n}{5}$, then $\psi_0-3=\tfrac{4(n-5)+3}{5}$ or $\psi_0-3=\frac{4(n-5)+5}{5}$, which contradicts the Theorem \ref{th2}; if $\frac{2n}{3}\le\psi_0\le\frac{4n-4}{5}$, then
\begin{align}
\Phi(T)&=2\cdot\Phi(T^{(1)})+\Phi(T^{(2)})\notag\\
       &\le 2\lambda^{4(n-5)-5(\psi-3)+2}+\lambda^{4(n-3)-5(\psi-2)+2}   \notag\\
         &=(\tfrac{2}{\lambda^5}+\tfrac{1}{\lambda^2})\lambda^{4n-5\psi+2}      \notag\\
         &<\lambda^{4n-5\psi+2},\notag
\end{align}
where we use $\lambda^3+2<\lambda^5$. This contradiction completes the proof of the case.\par
%Next we assume $\psi_0$ is odd. If $\psi_0=\frac{4n-1}{5}$ or $\frac{4n+1}{5}$, then $\psi_0-3=\frac{4(n-5)+4}{5}$ or $\psi_0-3=\frac{4(n-5)+6}{5}$, which contradicts the Theorem \ref{th2}; if $\frac{2n+1}{3}\le\psi_0\le\frac{4n-3}{5}$, then
%\begin{align}
%\Phi(T)&=2\Phi(T')+\Phi(T'')\notag\\
%       &=2f(\tfrac{4(n-5)-5(\psi_0-3)+2}{2})+g(\tfrac{4(n-3)-5(\psi_0-2)+1}{2})\notag\\
%       &=2f(\tfrac{4n-5\psi_0-3}{2})+g(\tfrac{4n-5\psi_0-1}{2})\notag\\
%       &\le2g(\tfrac{4n-5\psi_0-3}{2})+g(\tfrac{4n-5\psi_0-1}{2})\notag\\
%       &\le g(\tfrac{4n-5\psi_0+1}{2})\notag
%\end{align}
%where we use $f(n)\le g(n)$, which contradicts the choice of $T_0$.\par
\textbf{Case 2.} $\psi(T^{(1)})=\psi-4$.\par

In this case, there is at least one maximum dissociation set $F$ in $T^{(1)}$ such that $y\notin F$. On the other hand, every maximum dissociation set in $T$ contains $u$, $v$, $x$ and $w'$. So we have $\Phi(T)=\Phi_{\overline{y}}(T^{(1)})\le\Phi(T^{(1)})<\lambda^{4n-5\psi+2}$, which contradicts the choice of $T$.

We complete the proof of the claim.\qed
%
%First we assume $\psi_0$ is even. If $\frac{2n+2}{3}\le\psi_0\le\frac{4n}{5}$, then $\Phi(T)\le\Phi(T')\le f(\frac{4(n-5)-5(\psi_0-4)+2}{2})=f(\frac{4n-5\psi_0+2}{2})$, which contradicts the choice of $T_0$; if $\psi_0=\frac{2n}{3}$, then $\psi_0-4=\frac{2(n-5)-2}{3}$, which contradicts the Theorem \ref{th1}.\par
%
%Next we assume $\psi_0$ is odd. If $\psi_0=\frac{4n+1}{5}$, then $\psi_0-4=\frac{4(n-5)+1}{5}$. By the choice of $n$, the tree $T'$ has exactly one maximum dissociation set $F$ and for all vertex $v\in V(T')$ with $d_{T'}(v)\le2$ we have $v\in F$, which contradicts the condition of the case 2. If $\frac{2n+1}{3}\le\psi_0\le\frac{4n-1}{5}$, then $\Phi(T)\le\Phi(T')\le g(\frac{4(n-5)-5(\psi_0-4)+1}{2})=g(\frac{4n-5\psi_0+1}{2})$, which contradicts the choice of $T_0$; if $\psi_0=\frac{2n+1}{3}$, then by $\psi_0-4=\frac{2(n-5)-1}{3}$, which contradicts the Theorem \ref{th1}.\par\qed

\begin{cl}\label{cl3.6}
$x$ has no child that is a support vertex.
\end{cl}
\noindent{\bf Proof of Claim \ref{cl3.6}.}
Suppose, for a contradiction, that $x$ has a child $w'$ that is a support vertex. By Claim \ref{cl3.2} and \ref{cl3.3}, the vertex $w'$ has exactly one child $v'$ that is a leaf. Let $T^{(1)}=T-\{u,v\}$. Because it is easy to see that there exists at least one maximum dissociation set $F$ of $T$ such that $v\notin F$, $T^{(1)}$ has order $n-2$ and dissociation number $\psi-1$.
Let $T^{(2)}=T-\{u,v,w\}$, then $T^{(2)}$ has order $n-3$ and dissociation number $\psi-2$ and let $T^{(3)}=T-V(T_{x})$, then $T^{(3)}$ has order $n-6$ and dissociation number $\psi-4$.

%
%In this condition, there are $\Phi(T')$ maximum dissociation sets in $T_0$ that contain $u$ but not contain $v$; there are $\Phi(T''')$ maximum dissociation sets in $T_0$ that contain $v$ but not contain $u$ and there are $\Phi(T'')$ maximum dissociation sets in $T_0$ that contain $u$ and $v$. Since all maximum dissociation set in $T_0$ are one of these types, we have

A set $F$ is a maximum dissociation set of $T$ if and only if
\begin{itemize}
\item either $F=F'\cup \{u\}$, where $F'$ is a maximum dissociation set of cardinality $\psi-1$ of $T^{(1)}$,
\item or $F=F'\cup\{u,v\}$, where $F'$ is a maximum dissociation set of cardinality $\psi-2$ of $T^{(2)}$,
\item or $F=F'\cup\{v,w,v',w'\}$, where $F'$ is a maximum dissociation set of cardinality $\psi-4$ of $T^{(3)}$,
\end{itemize}
By the choice of $n$, this implies
\begin{align}
\Phi(T)&=\Phi(T^{(1)})+\Phi(T^{(2)})+\Phi(T^{(3)})\notag\\
          &\le \lambda^{4(n-2)-5(\psi-1)+2}+\lambda^{4(n-3)-5(\psi-2)+2}+\lambda^{4(n-6)-5(\psi-4)+2}  \notag\\
         &=(\tfrac{1}{\lambda^3}+\tfrac{1}{\lambda^2}+\tfrac{1}{\lambda^4})\lambda^{4n-5\psi+2}      \notag\\
         &=\lambda^{4n-5\psi+2},\notag
\end{align}
where we use $\lambda^2+\lambda+1=\lambda^4$. This contradiction completes the proof of the claim. \qed
%Furthermore, $\psi_0=\frac{2n}{3}$ implies that $\psi_0-2=\frac{2(n-3)}{3}$ and $\psi_0-4=\frac{2(n-6)}{3}$. If $\Phi(T_0)=f(\frac{4n-5\psi_0+2}{2})$ with $\psi_0=\frac{2n}{3}$, then $\Phi(T'')=f(\frac{4n-5\psi_0}{2})$ and $\Phi(T''')=f(\frac{4n-5\psi_0-2}{2})$, by the choice of $n$, which implies $T''=T(n-\psi_0-1)$ and $T'''=T(n-\psi_0-2)$, and hence $T_0=T(n-\psi_0)$, which contradicts the choice of $T_0$.\par
%
%Next we assume $\psi_0$ is odd. If $\psi_0=\frac{4n+1}{5}$ or $\psi_0=\frac{4n-1}{5}$, then $\psi_0-4=\frac{4(n-6)+5}{5}$ or $\psi_0-4=\frac{4(n-6)+3}{5}$, which contradicts the Theorem \ref{th2}; if $\frac{2n+1}{3}\le\psi_0\le\frac{4n-3}{5}$, by the choice of $n$, we obtain
%\begin{align}
%\Phi(T_0)&=\Phi(T')+\Phi(T'')+\Phi(T''')\notag\\
%         &\le f(\tfrac{4(n-2)-5(\psi_0-1)+2}{2})+g(\tfrac{4(n-3)-5(\psi_0-2)+1}{2})+g(\tfrac{4(n-6)-5(\psi_0-4)+1}{2})\notag\\
%         &=f(\tfrac{4n-5\psi_0-1}{2})+g(\tfrac{4n-5\psi_0-1}{2})+g(\tfrac{4n-5\psi_0-3}{2})\notag\\
%         &=g(\tfrac{4n-5\psi_0-1}{2})+2g(\tfrac{4n-5\psi_0-3}{2})+g(\tfrac{4n-5\psi_0-5}{2})\notag\\
%         &=g(\tfrac{4n-5\psi_0+1}{2})\notag
%\end{align}
%where we use $f(n)=g(n-1)+g(n-2)$, which contradicts the choice of $T_0$.\qed

\begin{cl}\label{cl3.7}
$d_{T}(x)=2$
\end{cl}
\noindent{\bf Proof of Claim \ref{cl3.7}.}
Suppose, for a contradiction, that $x$ has a child $w'$ distinct from $w$. By Claim \ref{cl3.5} and \ref{cl3.6}, the vertex $w'$ has a child $v'$ that has a child $u'$. By Claim \ref{cl3.2} and \ref{cl3.3}, $d_{T}(w')=d_{T}(v')=2$. Let $T^{(1)}=T-V(T_{x})$, then $T^{(1)}$ has order $n-7$ and dissociation number $\psi-4$ or $\psi-5$. Now, we consider the following two cases.

\textbf{Case 1.} $\psi(T^{(1)})=\psi-4$.\par
Let $T^{(2)}=T-(V(T_{x})\setminus\{x\})$, then $T^{(2)}$ has order $n-6$ and dissociation number $\psi-4$. In this case, $\psi(T^{(2)})=\psi(T^{(1)})$, by Lemma \ref{lem3.1}, we have
\begin{align}
\Phi_{\overline{x}}(T^{(2)})&\le\frac{1}{3}\Phi(T^{(2)})\label{31},\\
2\cdot\Phi_{\overline{x}}(T^{(2)})+\Phi_x^0(T^{(2)})&\le\Phi(T^{(2)})\label{32}.
\end{align}\par
A maximum dissociation set in $T^{(2)}$ not containing $x$ can be extended in nine ways to a maximum dissociation set in $T$, a maximum dissociation set $F$ in $T^{(2)}$ with $x\in F$ and $d_{T^{(2)}[F]}(x)=0$ can be extended in three ways to a maximum dissociation set in $T$, and a maximum dissociation set $F$ in $T^{(2)}$ with $x\in F$ and $d_{T^{(2)}[F]}(x)=1$ can only be extended in a unique way to a maximum dissociation set in $T$. Since all maximum dissociation sets of $T$ are of such forms, we obtain
\begin{align}
\Phi(T)&=9\cdot\Phi_{\overline{x}}(T^{(2)})+3\cdot \Phi_x^0(T^{(2)})+\Phi_x^1(T^{(2)}) \notag \\
              &=\Phi(T^{(2)})+8\cdot \Phi_{\overline{x}}(T^{(2)})+2\cdot \Phi_x^0(T^{(2)}). \notag
\end{align}
By inequalities (\ref{31}) and (\ref{32}), consider the following linear programming:
\begin{align}
\max.\ &8\cdot \Phi_{\overline{x}}(T^{(2)})+2\cdot \Phi_x^0(T^{(2)}) \notag \\
\text{s. t.} \ \ &\Phi_{\overline{x}}(T^{(2)}) \le \tfrac{1}{3}\Phi(T^{(2)}) \notag \\
&2\Phi_{\overline{x}}(T^{(2)})+\Phi_x^0(T^{(2)})\le \Phi(T^{(2)}) \notag \\
&\Phi_{\overline{x}}(T^{(2)})\ge 0,\ \ \Phi_x^0(T^{(2)})\ge0. \notag
\end{align}
The linear programming has an unique optimal solution
$$(\Phi_{\overline{x}}(T^{(2)}),\Phi_x^0(T^{(2)}))=(\tfrac{1}{3}\Phi(T^{(2)}),\tfrac{1}{3}\Phi(T^{(2)})).$$
Thus,
\begin{align}
\Phi(T_0)&\le \Phi(T^{(2)})+\tfrac{8}{3}\cdot \Phi(T^{(2)})+\tfrac{2}{3}\cdot \Phi(T^{(2)})\\
         &=\tfrac{13}{3}\cdot\Phi(T^{(2)})\notag\\
         &\le\tfrac{13}{3} \lambda^{4(n-6)-5(\psi-4)+2}\notag\\
         &<\lambda^{4n-5\psi+2},\notag
\end{align}
where we use $\tfrac{13}{3}<\lambda^4$. This contradiction completes the proof of the case. \par
%First, we assume $\psi_0$ is even. If $\frac{4n-4}{5}\le\psi_0\le\frac{4n}{5}$, then $\psi_0-4>\frac{4(n-7)+2}{5}$, which contradicts the Theorem \ref{th2}; if $\frac{2n}{3}\le\psi_0\le\frac{4n-6}{5}$, then by the choice of $n$, we have
%\[\Phi(T_0)\le\frac{13}{3}\Phi(T'')\le \frac{13}{3}f(\frac{4(n-6)-5(\psi_0-4)+2}{2})=\frac{13}{3}f(\frac{4n-5\psi_0-2}{2})\]
%By an inductive argument on $n$, for $n\ge4$, we have $f(n)\ge \frac{13}{3}f(n-2)$ with equality if and only if $n=4$. Specially, if $\psi_0=\frac{2n}{3}$ and $\frac{4n-5\psi_0+2}{2}=4$, then $n=9$ and $\psi_0=6$, which implies $T'=K_2$, hence $T_0=T(3)$, which contradicts the choice of $T_0$.\par
%
%Next we assume $\psi_0$ is odd. If $\frac{4n-5}{5}\le\psi_0\le\frac{4n+1}{5}$, then $\psi_0-4>\frac{4(n-7)+2}{5}$, which contradicts the Theorem \ref{th2}; if $\frac{2n+1}{3}\le\psi_0\le\frac{4n-7}{5}$, then by the choice of $n$, we have
%\[\Phi(T_0)\le\frac{13}{3}\Phi(T'')\le \frac{13}{3}g(\frac{4(n-6)-5(\psi_0-4)+1}{2})=\frac{13}{3}g(\frac{4n-5\psi_0-3}{2})\]
%By an inductive argument on $n$, for $n\ge4$, we have $g(n)> \frac{13}{3}g(n-2)$, which contradicts the choice of $T_0$.\par

\textbf{Case 2.} $\psi(T^{(1)})=\psi-5$.\par
In this case, there is at least one maximum dissociation set $F$ in $T^{(1)}$ such that either $y\notin F$ or $y\in F$ and $d_{T^{(1)}[F]}(y)=0$. Furthermore, a maximum dissociation set of $T^{(1)}$ can be extended at most three ways to a maximum dissociation set of $T$. Thus,
\begin{align}
\Phi(T)&\le 3\cdot\Phi(T^{(1)})\notag\\
         &\le3\lambda^{4(n-7)-5(\psi-5)+2}\notag\\
         &<\lambda^{4n-5\psi+2}\notag,
\end{align}
where we use $3<\lambda^3$. This contradiction completes the proof of the case.

We complete the proof of the claim.\qed

%First, we assume $\psi_0$ is even. If $\psi_0=\frac{4n}{5}$, then $\psi_0-5=\frac{4(n-7)+3}{5}$, which contradicts the Theorem \ref{th2}; if $\frac{2n+2}{3}\le\psi_0\le\frac{4n-2}{5}$, by the choice of $n$, we have $\Phi(T_0)\le3g(\frac{4(n-7)-5(\psi_0-5)+1}{2})=3g(\frac{4n-5\psi_0-2}{2})$. By an inductive argument on $n$, for $n\ge2$ we have $f(n)\ge 3g(n-2)$, which contradicts the choice of $T_0$; if $\psi_0=\frac{2n}{3}$, then $\psi_0-5=\frac{2(n-7)-1}{3}$, which contradicts the Theorem \ref{th1}.\par
%
%Next, we assume $\psi_0$ is odd. If $\psi_0=\frac{4n+1}{5}$, then $\psi_0-5=\frac{4(n-7)+4}{5}$, which contradicts the Theorem \ref{th2}; if $\psi_0=\frac{4n-1}{5}$, then $\psi_0-5=\frac{4(n-7)+2}{5}$. By the Theorem \ref{th2}, we know that the tree $T'$ has exactly one maximum dissociation set $F$ and all vertex $v$ in $T'$ with $d_{T'}(v)\le 2$ must be contained in $F$ with $d_{T'[F]}(y)=1$, which contradicts the condition of case 2; if $\frac{2n+1}{3}\le\psi_0\le\frac{4n-3}{5}$, by the choice of $n$, we have $\Phi(T_0)\le3f(\frac{4(n-7)-5(\psi_0-5)+2}{2})=3f(\frac{4n-5\psi_0-1}{2})$. By an inductive argument on $n$, for $n\ge2$, we have $g(n)> 3f(n-1)$, which contradicts the choice of $T_0$.\qed

We are now in a position to derive a final contradiction. By the above claims, $d_{T}(v)=d_{T}(w)=d_{T}(x)=2$. Let $T^{(1)}=T-V(T_{x})$, then $T^{(1)}$ has order $n-4$ and dissociation number $\psi-2$ or $\psi-3$. We consider the following two cases.\par

\textbf{Case 1.} $\psi(T^{(1)})=\psi-2$.\par
In this case, for every maximum dissociation set $F$ of $T^{(1)}$, we have $y\in F$ and $d_{T^{(1)}[F]}(y)=1$. Let $T^{(2)}=T-\{u,v\}$, then $T^{(2)}$ has order $n-2$ and dissociation number $\psi-1$; let $T^{(3)}=T-\{u,v,w\}$, then $T^{(3)}$ has order $n-3$ and dissociation number $\psi-2$.
A set $F$ is a maximum dissociation set of $T$ if and only if
\begin{itemize}
\item either $F=F'\cup\{v,w\}$, where $F'$ is a maximum dissociation set of cardinality $\psi-2$ of $T^{(1)}$,
\item or $F=F'\cup\{u\}$, where $F'$ is a maximum dissociation set of cardinality $\psi-1$ of $T^{(2)}$,
\item or $F=F'\cup\{u,v\}$, where $F'$ is a maximum dissociation set of cardinality $\psi-2$ of $T^{(3)}$.
\end{itemize}
%There are $\Phi(T'')$ maximum dissociation sets in $T_0$ that contain $u$ but not contain $v$; there are $\Phi(T')$ maximum dissociation sets in $T_0$ that contain $v$ but not contain $u$ and there are $\Phi(T''')$ maximum dissociation sets in $T_0$ that contain $u$ and $v$. Since all maximum dissociation set of $T_{0}$ are of such forms, we obtain
By the choice of $n$, this implies
\begin{align}
\Phi(T)&=\Phi(T^{(1)})+\Phi(T^{(2)})+\Phi(T^{(3)})\notag\\
         &\le \lambda^{4(n-4)-5(\psi-2)+2}+\lambda^{4(n-2)-5(\psi-1)+2}+\lambda^{4(n-3)-5(\psi-2)+2}  \notag\\
         &=(\tfrac{1}{\lambda^6}+\tfrac{1}{\lambda^3}+\tfrac{1}{\lambda^2})\lambda^{4n-5\psi+2}      \notag\\
         &<\lambda^{4n-5\psi+2}\notag,
\end{align}
where we use $\lambda^4+\lambda^3+1<\lambda^6$. This contradiction completes the proof of the case.\par

\textbf{Case 2.} $\psi(T^{(1)})=\psi-3$.\par
In this case, there is at least one maximum dissociation set $F$ in $T^{(1)}$ such that either $y\notin F$ or $y\in F$ and $d_{T^{(1)}[F]}(y)=0$. Let $T'=T-\{uv\}+\{xu\}$, then it is easy to see that $T'$ has order $n$ and dissociation number $\psi$. Next, we will prove that $\Phi(T)\leq\Phi(T')$.

A set $F$ is a maximum dissociation set of $T$ if and only if
\begin{itemize}
\item either $F\in\{F'\cup\{x,v,u\},F'\cup\{x,w,u\}\}$, where $F'$ is a maximum dissociation set of cardinality $\psi-3$ of $T^{(1)}$ such that $y\notin F'$,
\item or $F=F'\cup\{x,v,u\}$, where $F'$ is a maximum dissociation set of cardinality $\psi-3$ of $T^{(1)}$ such that $y\in F'$ and $d_{T^{(1)}[F']}(y)=0$.
\end{itemize}
Thus, we have $\Phi(T)=2\cdot\Phi_{\overline{y}}(T^{(1)})+\Phi^0_y(T^{(1)})$.

Similarly, a set $F$ is a maximum dissociation set of $T'$ if and only if
\begin{itemize}
\item either $F\in\{F'\cup\{x,v,u\},F'\cup\{w,v,u\}\}$, where $F'$ is a maximum dissociation set of cardinality $\psi-3$ of $T^{(1)}$ such that $y\notin F'$,
\item or $F=F'\cup\{w,v,u\}$, where $F'$ is a maximum dissociation set of cardinality $\psi-3$ of $T^{(1)}$ such that $y\in F'$ and $d_{T^{(1)}[F']}(y)=0$,
\item or $F=F'\cup\{w,v,u\}$, where $F'$ is a maximum dissociation set of cardinality $\psi-3$ of $T^{(1)}$ such that $y\in F'$ and $d_{T^{(1)}[F']}(y)=1$.
\end{itemize}
Thus, we have $\Phi(T')=2\cdot\Phi_{\overline{y}}(T^{(1)})+\Phi^0_y(T^{(1)})+\Phi^1_y(T^{(1)})$. It follows that $\Phi(T)\leq\Phi(T')$. On the other hand, by Claim \ref{cl3.3},
\begin{align}
\Phi(T)\le\Phi(T')\le \lambda^{4n-5\psi+2}.\notag
\end{align}
%In the subcubic tree $T_0$, a maximum dissociation set in $T'$ not containing $y$ can be extended in two ways to a maximum dissociation set in $T_0$ and every maximum dissociation set in $T_0$ not containing $y$ is of that form; a maximum dissociation set in $T'$ containing $y$ with $d_{T'[F]}(y)=0$ can be extended in an unique way to a maximum dissociation set in $T_0$ and every maximum dissociation set in $T_0$ containing $y$ is of that form.\par
%
%Similarly, in the subcubic tree $T_1$, a maximum dissociation set in $T'$ not containing $y$ can be extended in two ways to a maximum dissociation set in $T_1$ and every maximum dissociation set in $T_1$ not containing $y$ is of that form; a maximum dissociation set in $T'$ that contains $y$ with $d_{T'[F]}(y)=0$ can be extended in an unique way to a maximum dissociation set in $T_1$ and a maximum dissociation set in $T'$ that contains $y$ with $d_{T'[F]}(y)=1$ can be extended an unique way to a maximum dissociation set in $T_1$. Every maximum dissociation set in $T_1$ containing $y$ is of that form.
 %% Arguing as in the proof of the claim \ref{cl3.3}, we have
This final contradiction completes the proof.\qed

\begin{coro}
Let $T$ be a subcubic tree of order $n$ and dissociation number $\psi$. If $\psi=\tfrac{4n+2}{5}$, then $T$ has exactly one maximum dissociation set, i.e., $\Phi(T)=1$.
\end{coro}

\pf The result can be easily obtained from Theorem \ref{th3.2}. In fact, if $\psi=\tfrac{4n+2}{5}$, by Theorem \ref{th2.3}, the structure of $T$ is known and it can been seen that $\Phi(T)=1$. \qed

\begin{coro}
Let $T$ be a subcubic tree of order $n$. Then
$$\Phi(T)\leq 1.29^{n+1}.$$
\end{coro}

\pf By Theorem \ref{th2.4} and \ref{th3.2}, we have

$$\Phi(T)\leq \lambda^{4n-5\cdot \tfrac{2n}{3}+2}=\lambda^{\tfrac{2n+2}{3}}\leq1.29^{n+1}.$$ \qed

\begin{coro}
Let $T$ be a subcubic tree of dissociation number $\psi$. Then
$$\Phi(T)\leq 1.466^{\psi+2}.$$
\end{coro}

\pf By Theorem \ref{th2.4} and \ref{th3.2}, we have

$$\Phi(T)\leq \lambda^{4\cdot\tfrac{3\psi}{2}-5\psi+2}=\lambda^{\psi+2}=1.466^{\psi+2}.$$ \qed

\section{Further remark}

In fact, with some tedious calculation we can obtain a more accurate value on the maximum number of maximum dissociation sets of a subcubic tree of order $n$ and dissociation number $\psi$. Furthermore, the extremal subcubic trees achieving this maximum are also found.

Let $f(n)$ and $g(n)$ be two sequences, where $f(0)=1$, $f(1)=1$, $f(2)=3$ and $f(n)=f(n-1)+2f(n-2)+f(n-3)$ when $n\ge3$ and $g(0)=1$, $g(1)=2$, $g(2)=4$ and $g(n)=g(n-1)+2g(n-2)+g(n-3)$ when $n\ge3$. The equation $\alpha^3-\alpha^2-2\alpha-1=0$ mentioned in the proof of Theorem \ref{th3.2} is exactly the characteristic equation of $f(n)$ and $g(n)$.

For a positive integer $\ell$,

(1) let $T_{\ell}$ arise by attaching a pendant edge to every vertex of a path of order $\ell$, then it can be seen that $V(T_{\ell})=3\ell$, $\psi(T_{\ell})=2\ell$ and $\Phi(T_{\ell})=f(\ell+1)$;

(2) let $T^{(1)}_{\ell}$ arise by attaching a pendant edge to every vertex of a path of order $\ell+1$ except exactly one of its two endvertices, then it can be seen that $V(T^{(1)}_{\ell})=3\ell+1$, $\psi(T^{(1)}_{\ell})=2\ell+1$ and $\Phi(T^{(1)}_{\ell})=g(\ell)$;

(3) and let $T^{(2)}_{\ell}$ arise by attaching a pendant edge to every vertex of a path of order $\ell+2$ except its two endvertices, then it can be seen that $V(T^{(2)}_{\ell})=3\ell+2$, $\psi(^{(2)}_{\ell})=2\ell+2$ and $\Phi(T^{(2)}_{\ell})=f(\ell)$.

For two positive integers $n$ and $\psi$ with $\frac{2n}{3}\le\psi\le\frac{4n+2}{5}$, let $T(n,\psi)$ denote the extremal subcubic tree having the maximum number of maximum dissociation sets among all subcubic trees of order $n$ and dissociation number $\psi$. Using the proof method in the proof of Theorem \ref{th3.2} and a more complex analysis, the following result can be obtained.
\begin{theo}
If $T$ is a subcubic tree of order $n$ and dissociation number $\psi$, then,
\begin{center}
$\Phi(T)\le\left \{
  \begin{array}{ll}
   f(\frac{4n-5\psi+2}{2})&,  \hbox{if $\psi$ is even;}\vspace{8pt}\\
   g(\frac{4n-5\psi+1}{2})&, \hbox{if $\psi$ is odd.}
\end{array}
\right.$
\end{center}
Furthermore,

(1) if $\psi=\frac{2n}{3}$, then $T(n,\psi)=T_{\psi/2}$;

(2) if $\frac{2n}{3}<\psi<\frac{4n}{5}$ and $\psi$ is even, then $T(n,\psi)$ is obtained either from $T_{\frac{4n-5\psi}{2}}$ by attaching $\frac{3\psi-2n}{2}$ times $P_5s$ to some vertex of degree less than 3, or from $T^{(2)}_{\frac{4n-5\psi+2}{2}}$ by attaching $\frac{3\psi-2n-2}{2}$ times $P_5s$ to some vertex of degree less than 3;

(3) if $\frac{2n}{3}<\psi<\frac{4n}{5}$ and $\psi$ is odd, then $T(n,\psi)$ is obtained from $T^{(1)}_{\frac{4n-5\psi+1}{2}}$ by attaching $\frac{3\psi-2n-1}{2}$ times $P_5s$ to some vertex of degree less than 3;

(4) if $\psi=\frac{4n}{5}$, then $T(n,\psi)$ is obtained from $P_5$ or $T^{(2)}_{1}$  by attaching $\frac{3\psi-2n-1}{2}$ times $P_5s$ to some vertex of degree less than 3;

(5) if $\psi=\frac{4n+1}{5}$, then $T(n,\psi)$ is obtained from $K_1$ by attaching $n-\psi$ times $P_5s$ to some vertex of degree less than 3;

(6) if $\psi=\frac{4n+2}{5}$, then $T(n,\psi)$ is obtained from $K_2$ by attaching $n-\psi$ times $P_5s$ to some vertex of degree less than 3.

%It is easy to check that a tree $T\in\mathcal{T}(n,\psi)$ has order $n$ and dissociation number $\psi$. Since attaching $P_5$ to a tree do not change the number of maximum dissociation set of the tree, we can obtain that for a tree $T\in\mathcal{T}(n,\theta)$, $\Phi(T)=f(\frac{4n-5\psi+2}{2})$ and $\Phi(T)=g(\frac{4n-5\psi+1}{2})$, where $\psi$ is even and odd respectively.
\end{theo}

\section*{Acknowledgment}
This work was supported by National Key Research and Development Program of China (2019YFC1906102), National Key Technology Research and Development Program of the Ministry of Science and Technology of China (No. 2015BAK39B00).

\bibliographystyle{unsrt}
%\bibliography{acl}

\end{document}